\documentclass[12pt]{article}
\usepackage{geometry,amssymb,amsthm,amsmath,
graphics,enumerate}

\let\oldthebibliography\thebibliography
\renewcommand{\thebibliography}[1]{%
  \oldthebibliography{#1}%
  \setlength{\itemsep}{0pt}%
  \setlength{\parskip}{0pt}%
}

\usepackage{xcolor}
\usepackage{times}
\usepackage{amsfonts}
\usepackage{amscd}
\usepackage{url}
\usepackage{showlabels}

  \usepackage[normalem]{ulem}
\usepackage{soul}
\usepackage{hyperref}

\newcommand{\sidebar}[1]{\vskip10pt\noindent
 \hskip.70truein\vrule width2.0pt\hskip.5em
 \vbox{\hsize= 4truein\noindent\footnotesize\relax #1 }\vskip10pt\noindent}

\newbox\smilebox
\newbox\anchorbox
\newbox\noanchorbox
\newbox\tempbox

\setbox\smilebox=\hbox{$\smile$}

\def\anchor{\hbox{\vtop{
           \hbox to \wd\smilebox{\hfil\vrule width.4pt height7pt depth1pt\hfil}
           \vskip  -11.5truept
           \hbox to \wd\smilebox{\hfil$\smile$\hfil}}}}
\setbox\anchorbox=\anchor
\def\noanchor{\hbox{\vtop{
           \hbox to \wd\anchorbox{\hfil\anchor\hfil}
           \vskip -14truept
           \hbox to \wd\anchorbox{\hfil/\hfil}}}}
\setbox\noanchorbox=\noanchor

\def\fg#1#2#3{\setbox\tempbox=\hbox{$\scriptstyle{#2}$}
\ifnum\wd\anchorbox>\wd\tempbox\dimen255=\wd\anchorbox
\else\dimen255=\wd\tempbox\fi
{#1\,\vtop{\hbox to \dimen255{\hfil\anchor\hfil}
           \vskip -6truept
           \hbox to \dimen255{\hfil$\scriptstyle{#2}$\hfil}}
           \,#3}}

\def\nfg#1#2#3{\setbox\tempbox=\hbox{$\scriptstyle{#2}$}
\ifnum\wd\noanchorbox>\wd\tempbox\dimen255=\wd\noanchorbox
\else\dimen255=\wd\tempbox\fi
{#1\,\vtop{\hbox to \dimen255{\hfil\noanchor\hfil}
           \vskip -6truept
           \hbox to \dimen255{\hfil$\scriptstyle{#2}$\hfil}}
           \,#3}}

\setbox1=\hbox{$\bot$}

\def\north#1#2{#1\,
\hbox{$\bot$\llap {\hbox to\wd1 {\hfil $/$\hfil}}}
\,#2}

\def\nao#1#2#3{#1\  \hbox{\vtop{
\baselineskip=4pt
\hbox{$\bot$\llap {\hbox to\wd1 {\hfil $/$\hfil}}
\hskip .05em \llap{\hbox{$^{\scriptscriptstyle{a}}$}}}\hbox{$\scriptstyle
{#2}$}}}\, #3}

\def\bbar{\overline{b}}

\def\dbar{\overline{d}}

\def\hbar{\overline{h}}

\def\vbar{\overline{v}}

\def\Qbar{\bf Q}
\def\Rbar{\bf R}

\def\C{{\mathfrak  C}}

\def\FF{{\bf F}}
\def\C{{\cal  C}}

\def\FF{{\bf F}}

\def\LL{{\mathcal  L}}

\def\acl{{\rm acl}}

\def\Fa0{{\FF^a_{\aleph_0}}}

\def\<{\langle}

\def\>{\rangle}

\def\Lscr{\mathcal L}

\def\phi{\varphi}
\def\PA{PA}

\def\Const{Const}
\def\C_{C_}

\def\DN_{DN}

\makeatother
\newtheorem{theorem}{Theorem}[section]

\newtheorem{definition}[theorem]{Definition}

\newtheorem{remark}[theorem]{Remark}

\newtheorem{corollary}[theorem]{Corollary}

\newtheorem{procedure}[theorem]{Procedure}

\begin{document}

\author{John T. Baldwin, Constantin C. Br\^incu\c s}

\title{Carnapian Frameworks and Categoricity of Arithmetic via Inferential $\omega$-logics}
\maketitle

\begin{abstract} 
 We provided in \cite{BaldwinBrincusI} extensions of first order logic by  modified inferential definitions of the classical $\omega$-rule
  in $1$ or $2$ sorts. These logics are categorical in the inferential sense. Arithmetic has a unique countable model in each case, e.g. first order PA is categorical in our first logic. The 2-sorted case
  interprets $L_{\omega_1,\omega}$. In this paper, we discuss two philosophical problems raised by Button and Walsh \cite{ButtonWalshbook} concerting  the identification of a unique isomorphism class. First, we argue that the doxological challenge (on referential determinacy) gets a clear answer if placed in an appropriate (Carnapian) linguistic framework and   is meaningless otherwise. To clarify this approach,
  we address Button-Walsh's dismissal of concepts-modelism
  by developing the notion of {\em cognitive modelism}, according to which classical mathematics is a complex process of constructing and developing a distinctive class of concepts.
    Second, we argue that the inferential $\omega$-logics, that are much weaker than second order logic, do not appeal to the arithmetical concepts that the categoricity theorems 
    proved within these logics aim to secure.

  

\end{abstract}

\section{ Introduction}\label{intro}
Our \cite{BaldwinBrincusI}
provides  modified inferential definitions of the classical $\omega$-rule that extend  $L_{\omega,\omega}$. 
The key features of these rules are that they allow (putative) non-standard parameters to appear in  the formula in the hypothesis  and yield a categorical interpretation of the universal quantifier.
Extending first order logic by{\em Our Inferential $\omega$-rule} (the $I-\omega$-rule) 
yields Peano arithmetic (even more, Robinson's system $\Qbar$) is categorical.
 Working in two-sorted $G-\omega$-logic ($G$ for generalized) we  proved that each 
$L_{\omega_1,\omega}(\tau)$-complete (Scott) sentence $\phi$ is `structurally equivalent with' (same $\tau$-models as) a pseudo-elementary class (the relativized $\tau$-reducts of an associated first-order theory ${\hat T}^\phi$ in the $G-\omega$-logic).  This reduces the philosophical issue of the existence of infinite formulas to that of the reliability of infinitary proofs in a first order context. Further,
 a Scott sentence is categorical if and only if the associated ${\hat T}^\phi$   is in the $G-\omega$-logic and such structures are chararacterized. 
We briefly introduce in Section \ref{Infomega} the main notions and results from \cite{BaldwinBrincusI} that are useful for the philosophical discussion in the sections below.

Based on these results, we discuss in this paper two philosophical problems raised by Button and Walsh \cite{ButtonWalshbook} in connection to the identification of a unique isomorphism class, i.e. the doxological challenge (Section \ref{Dochall}) and the circularity argument (Section \ref{AgCA}).
 The doxological challenge asks a  philosopher: `How can we pick out an isomorphism type'. They argue that various species of philosopher (modelists of different varieties) cannot answer this question.   The circularity argument asserts that any proof of a categoricity theorem  for arithmetic invokes mathematical concepts that were supposed to be justified by that very categoricity theorem. We argue (in Section \ref{carling}) that the purported difficulty of the doxological challenge is generated by the confusion between the notions of {\em structure for a vocabulary (signature)} and {\em model of a theory} and we clarify it by means of four (Carnapian) frameworks (i.e. structure framework, theoretical framework of set theory, metatheoretical framework of formalized theories in  systems of logic, and the philosophical framework). The challenge proves meaningful and easily answered in the second framework, but it lacks cognitive content, and it is thus meaningless in the first one. In Section~\ref{cogmodelism}, building on such as \cite{DeBenRos, Carey, Devlingene,LakNun} we develop {\em cognitive modelism} as a response to their dismissal of `concepts modelism' \cite[p 149]{ButtonWalshbook}. 
 We address the circularity argument in Section~\ref{AgCA} by showing that our inferential approach requires only the notions of finite and countably infinite (equicardinality) and so does not presuppose the notion of an $\omega$-sequence or arithmetical notions.

\section{Categoricity by inferential $\omega$-logic and  $L_{\omega_1,\omega}$}\label{Infomega}

In this section we state the major results of \cite{BaldwinBrincusI} whose philosophical consequences are analyzed in the bulk of this paper.

\begin{definition}[Vocabularies and structures]\label{str}

By a  vocabulary $\tau$ we mean a set of relation and function symbols with prescribed arity. A $0$-ary
function is called a constant. 
These are not logical symbols,
but vary with the topic being formalized. For any vocabulary $\tau$, a $\tau$-structure $M$ consists of a domain $M$  equipped with: for
each $n$-ary relation symbol $R$  in $\tau$ a subset $R^M \subseteq M^n $, the sequences that satisfy $R$,
and for each $n$-function symbol $f$ a function $f^M$ mapping $M^n$ into $M$ (see \cite[\S 1.1]{ButtonWalshbook}, \cite{Markerbook}).

\end{definition}

\begin{remark}[$\omega$-logic, $\omega$-rule, $\omega$-model and standard model] \label{omega-logics and rules}
{\rm In all cases we let $c$ range over a countable set $C$ of constants. 
\begin{enumerate}
\item $\omega$-logic and $\omega$-rule. 
There are both  model-theoretic   and   proof-theoretic approaches to $\omega$-logic.  Chang and Keisler (\cite[p.82-3]{ChangKeisler}) expound both variants. 
\begin{enumerate}
\item  
The more traditional proof-theoretic one  adds the `classical' $\omega$-rule to the rules of inference of {\em one-sorted} first-order logic:\footnote{In arithmetic contexts the constants are often called `standard numerals'. While this notion developed from the study of arithmetic, this generalization has no reliance on arithmetic.} 

$$\omega\text{-rule}:  \frac{\{\phi (c): c \in C\}}{(\forall x)\phi (x)}$$\label{classwrule}

 Precise inferential formulation of  this rule will be given in Definition ~\ref{ourinfrule}.

 \item The model-theoretic version, stemming from (\cite{Henkinomt, Orey, MorleyLeeds} as treated in  \cite[pp.28, 39]{Ebbing}, and \cite[pp. 153-155] {Shapiro}) is {\em two-sorted} with a designated predicate for `the natural numbers'.    In this case, an $\omega$-model $M$ is one where the predicate $N(M)$ consists only of a designated set of constants.  

  And the $\omega$-rule is:
  $$\text{Generalized}\hspace{.1in} \omega\text{-rule}:  \frac{\{\phi (c): c \in C\}}{(\forall x )(N(x) \rightarrow\phi(x))}.$$
We modify this  rule to the $G-\omega$-rule in  Definition~\ref{ourinfrule2}.

  \end{enumerate}

\item $\omega$-model: There are several meaning of this term; the most important distinction  here is whether the context is $1$-sorted or $2$-sorted.
\begin{enumerate}[(a)]
\item $1$-sorted:
\begin{enumerate}[(i)]
 \item {\em $\omega$-sequence  $\omega$-model}
 In the vocabulary $(0,S)$ the domain of the $\omega$-model
 is the set of iterations of applications of a $1-1$ function to $0$.

    \item {\em arithmetic $\omega$-model}
The concept arose in studying arithmetic and the classical notion is: $M$ is an $\omega$-model of $T$ if $T$ interprets arithmetic in the vocabulary $(+,\times,0,1)$ of arithmetic and 
the universe of $M$ is named by the numerals.

\item {\em first order $\omega$-model}
This is generalized by replacing `natural numbers' with `the denotations of a specified countable set of constants' and $T$ is any first order theory. 

\end{enumerate}
\item  $2$-sorted: \cite{Orey} introduced the $2$-sorted approach with a predicate $N$, a vocabulary containing countably many constants from $N$
and an $\omega$-model is one where these constants exhaust $N$. 

\end{enumerate}

     \item Standard model:  A {\em standard model} arises when there is an informal notion that has  widely accepted formal counterpart.  E.g.,  arithmetic,  consider a vocabulary $\tau$ that includes the following non-logical terms $\{0, S, +, \times\,, <\}$ with the following abbreviations: $1=S0, 2=SS0, \ldots$. A standard model of arithmetic is an $\omega$-model in which the domain is $N=\{0,1,2, 3,\ldots\}$, i.e. the model omits the set of formulas $\{x\neq 0, x\neq 1, x\neq 2, \ldots\}$. As \cite[p. 28, 5B]{ButtonWalshbook} point out, a confusion may arise when, influenced by Dedekind's axioms and the omnipotence of second order definability, one takes only $(\omega, 0,S)$ as the standard model of arithmetic.
       
\end{enumerate}}

\end{remark}

The crucial innovation for the one-sorted case in \cite{BaldwinBrincusI} is:

\begin{definition}{(Inferential $\omega$-rule)}\label{ourinfrule} Fix a vocabulary $\tau$ and use the notation of  Robinson-semantics (\cite[2.4.2]{BaldwinBrincusI}, \cite[pp.15-19]{ButtonWalshbook}). For  any countable set of constants $C\supseteq \Const(\tau)$ and any  well-formed formula $\phi(x,\dbar)$ in the vocabulary  $\tau$ expanded by the constants from $C- \Const(\tau)$:

\[I-\omega\text{-rule}: \frac{\bigwedge\{\phi(c ,\dbar): c \in \Const(\tau)\}} {(\forall x)\phi(x,\dbar)};\]

$$I-\forall E:       \frac{(\forall x) \phi(x,\dbar)}
        {\phi (c, \dbar), \text{ for each } c\in C}.$$
\end{definition}

And here is the important consequence.
\begin{theorem} \label{wrulebworks} 


If a first order theory $T$ in a vocabulary $\tau$ has a
model  
that embeds in every model of $T$ (is algebraically prime)
and with every element named by a $\tau$-constant, $T$  plus the $I-\omega$-rule (Definition~\ref{ourinfrule})
is categorical. 

Two examples are 
Robinson's $\Qbar$  and the theory of $(\omega,0,S)$
with axioms that $S$ is $1-1$ and every point is in the range except $0$.
\end{theorem}


 Without algebraically prime, the result is slightly weaker.

\begin{theorem}\label{allmax} Fix a vocabulary $\tau$ with $\aleph_0$ constants and suppose that for any model of $T$, the substructure consisting of the elements named by those constants satisfies $T$ and can be embedded in every
model of $T$. Then,  no model of $T^{\omega}$ has a proper extension.
\end{theorem}

We now move to the two-sorted case and describe the theory $T^{\hat\phi}$ that represents models of an $L_{\omega_1,\omega}$ sentence $\phi$.

\begin{definition}\label{Morcode} {Fix a Scott sentence\footnote{I.e. All models of $\phi$ are $\infty,\omega$-equivalent; equivalently $\phi$  is $\aleph_0$-categorical.}  $\phi$ in a vocabulary $\tau$ with a 
countable model $A$.  
We work in a two-sorted\footnote{Constants and variable will be restricted to specific sorts.} vocabulary $\sigma$ ($\sigma$ depends on $\tau$) with disjoint sorts $(N,V)$ where $N$ contains the set of images of  the constants $\Const(N)=_{\rm df} \langle c_{n,i}:i,n <\omega\rangle$ and $V$ consists of a
$\tau^\phi$-structure satisfying the theory $T^\phi$ constructed from $\phi$ by Chang\footnote{ Chang 
\cite[4.1.2]{BaldwinBrincusI} constructs from a $\tau$-Scott sentence $\phi$ a first order $\tau^\phi$ theory whose models reducted to $\tau$ are the models of $\phi$.}
Each $\tau^\phi$-constant becomes a $\sigma$-constant satisfying $V$.

Adapting \cite{MorleyLeeds}\footnote{
The constants $c_{n,i}$ satisfying $N(x)$ code all $T^{\phi}$ principal types over $\emptyset$.}, we construct a theory 
 $T^{\hat \phi}$ such that a  $G-\omega$-model
of $T^{\hat \phi}$ omits each non-principal type over $\emptyset$.  In particular, if $B\models T^{\hat \phi}$, the restriction of $B$ to $V(B)$
satisfies $T^{\phi}$. 
Extend the  vocabulary $\sigma$ of $T^{\hat \phi}$ to include,
for each   $n$, an $(n+1)$-ary-relation $R^n$ on $N\times V^n$  in 2-sorted G-$\omega$-logic.
 
The theory $T^{\hat \phi}$
is given by $T^{\phi}$ relativized to $V$ along with axioms saying the $c_{n,i}$ are distinct elements of $N$;
the following axioms on the $R^n$ ensure  that the elements of $N$ code the atomic $\tau^\phi$-types  over the empty set of the elements of $V$. 
  $$ (1)\  (\forall v_0) [R^n(c_{n,i} ,\vbar_0) \leftrightarrow \phi_{i}(\vbar_0)].$$ where
 $\phi_i(\vbar)$ generates the $i$th principal  $n$-type (in $\tau^\phi$) over $\emptyset$ for a given injective enumeration of those types\footnote{ $R^n(c_{n,i},\dbar)$ with $\dbar \in V$ means: $\dbar$ realizes the   $n-\tau^\phi$-type over the empty set indexed by $c_{n,i}$.}. 
    $$(2) \ (\forall \vbar_0) [V(\vbar_0) \rightarrow \exists v_1 [N(v_1)\wedge R^n(v_1,\vbar_0) ]$$
We call a model $B$ of $T^{\hat \phi}$ with $\Const(N)= \{c_{n,i}:i,n<\omega\}$ denoting the elements of $N(B)$ a $G-\omega$-model. 
}
\end{definition}

And now we describe the $\omega$-rule that makes that interpretation work. Again, the key is extending the constants allowed in the active formula.

\begin{definition}[Inferential $G-\omega$-rule]\label{ourinfrule2} 
$\sigma$ is a $2$-sorted vocabulary as in Definition~\ref{Morcode}; For  any countable set of constants $B\supseteq \Const(\sigma)$ and for any 
$\sigma(B)$-formula $\lambda(x,\vbar)$ and,  in particular, instances $\psi(x,\dbar)$ of a formula $\psi(x,\vbar)$  with $\lg(\vbar) =n$:

{\small$$
G\text{-}\omega\text{-rule}:\;
\frac{\bigwedge\{\psi(c_{n,i},\dbar): c_{n,i}\in \Const(N), n,i<\omega\}
\&
(\forall x)(\forall v_i)(\psi(x,\dbar) \rightarrow (N(x) \& \bigwedge_{i< n}V(v_i)))}
{(\forall x)(N(x) \rightarrow \psi(x,\dbar))}
$$}

$$G-\forall E:       \frac{(\forall x) \lambda(x,\bbar)}
{\lambda (e,\bbar), \text{ for each } e \in B}.$$

Note that there are constants in $B$ that are not in $N$, but they are not instances of the formula in the hypotheses of the $G-\omega$-rule.

\end{definition}
 \begin{theorem}\label{omegaruleatomic2} Fix a complete $L_{\omega_1,\omega}$-sentence $\phi$. 
With Definition~\ref{ourinfrule2}, if a model of $T^{\hat \phi}$ satisfies the $G-\omega$-rule 
then its  $\tau^\phi$-reduct   is an atomic model of $T^{ \phi}$. 
I.e., each  $p\in S(\emptyset)$ realized in $V$  is principal as a $\tau^\phi$-type.  Thus, the $\tau$-reducts of models of $T^{\hat \phi} + G-\omega-rule$ are exactly the models of $\phi$. That is, $\phi$ and $T^{\hat \phi} + G-\omega$-rule are structurally equivalent.
\end{theorem}

\begin{definition}\label{gendef} A structure $M$ is said to be generative if it is isomorphic to a proper substructure (equivalently extension)
of itself.
\end{definition}

\begin{corollary}\label{coupdegras} $\phi$ is categorical if and only if $T^{\hat \phi}$ is categorical in $G-\omega$-logic if and only if the countable model of $\phi$ is non-generative.
\end{corollary}

\section{Structures, Models, and the Modelists Doxological Challenge}\label{Dochall}

In proving the categoricity results in the inferential $\omega$-logics, we introduced in Definition ~\ref{str} a notion of structure that is identical {\em mutatis mutandis} to the one defined by \cite[Definition 1.2]{ButtonWalshbook}. 
As \cite[p. 8]{ButtonWalshbook} suggest,  writing `as usual in set theory', this notion of structure is given in (naive?) set theory. Our aim in proving categoricity for {\em theories} was to choose (or given) a vocabulary $\tau$ and a $\tau$-structure, $M$, find a logic $\Lscr$ and an $\Lscr(\tau)$-theory $T$ whose unique model up to isomorphism is $M$. We aimed to find substantially weaker choices of $\Lscr$ that can give categorical of appropriate countable structures.

Although they are situated in different frameworks (Cf. \S \ref{carling}), the notions of {\em structure} and {\em model of a theory} are not always clearly and coherently kept distinct, and this generates philosophical misunderstandings. For instance, even though they define the notion of $\tau$-structure both informally just before (\cite[Definition 1.2]{ButtonWalshbook}) and then semiformally  in set theory exactly as in a model theory text, \cite{ButtonWalshbook} formulate a problem ({\em The Modelist's Doxological Challenge}) that obscures the difference between the informal notion of structure and the formal notion of model\footnote{The notion of structure is well-defined with respect to a specific signature, but, unlike `model of a theory' it is not dependent on any formal system of logic. }, and likewise between the different frameworks to which these two notions belong. We consider the context of the  challenge and provide an answer to it.  For the clarification, we specify `frameworks' in \S \ref{carling} to distinguish between the identification and  uniqueness of structures in informal but formalizable mathematics (Framework 2 in \S 3) -- on the one side 
%
and categoricity (uniqueness up to isomorphism) of models of a formal theory --on the other side (framework 3 ).


\subsection{The Modelists and the Doxological Challenge}
\label{moddc}

\cite[ p. 143-4]{ButtonWalshbook} introduced the notion of modelism, as an attitude towards model-theory, with this claim as to how `structure' should be explicated:

\begin{quote}{\em The modelists manifesto.} Mathematicians frequently engage in structure-talk, and model theory 
provides precise tools for explicating the notion of {\em structure} they invoke. Specifically, mathematical structure should be understood in terms of {\em isomorphism}, in the model theorist’s sense.

\end{quote}

This commandment brings joy to the ears of a model-theorist, especially since both the notions of structure and isomorphism (\cite[Definition 2.2]{ButtonWalshbook}) are defined in naive set theory exactly as in  model-theory texts. However, it is then rather surprising to read:

\begin{quote}{\em The Modelist's Doxological Challenge.}  How can we pick out particular isomorphism types?
\cite[p 145]{ButtonWalshbook}
\end{quote}

 The challenge is raised for various kinds 
 of modelists that \cite[144-148]{ButtonWalshbook} consider. Each of these modelists accepts that mathematicians deal with structures\footnote{\cite{Reck} identifies five forms of structuralism: methodological, set-theoretic, ante rem, modal, and logical. This taxonomy is more fine-grained and partly orthogonal to ~\cite{ButtonWalshbook}'s classification. Methodological structuralism, according to ~\cite[p.371]{Reck}
 emphasizes a formal, axiomatic method of presentation of mathematics rather than metaphysical commitments. We can think of it as background for the distinction between set-theoretical and logical structuralism. Reck's set-theoretic structuralism mirrors what we have called the model-theoretic approaches, but he adds certain metaphysical and semantic theses. 
His distinction between set-theoretic and logical structuralism is significant here. Dedekind’s account exemplifies the latter: he takes $N$ as any representative of the isomorphism type and then `creates' the natural numbers, in a {\em logical} manner. The {\em constitutive properties} of these numbers consist only in those imposed by successor. While other members of the isomorphism type (e.g. the von Neumann ordinals) have other constitutive properties evident in the set theoretic construction of the structure.}. Blunt object-modelists claim that structure-talk is about proper (abstract) objects, while concepts-modelists maintain that this talk is about mathematical concepts. \cite[146]{ButtonWalshbook} argue that each  modelist fails to explain either i) how we refer to {\em bona fide} objects or ii) how we acquire and use concepts that are as fine-grained as isomorphism types, respectively.

 We take the doxological challenge to be only partly meaningful, since we consider it to be based on the confusion between internal and external questions that \cite{Carext}  warned us about long ago. Briefly, once we set up a certain theoretical framework, the question either is internal to that framework, and thus gets a precise answer within the framework, or it is external to it, and it is thus meaningless, i.e. it lacks cognitive content. 
 We expand this discussion of i) in Section \S~\ref{carling}. 
 \cite{DeBenRos}, based on evidence from Cognitive Science, advance convincing arguments towards answering ii) which we expand in Section
 \S~\ref{cogmodelism}.

Writing $M \approx N$ to indicate that there is a bijection between $M$ and $N$ that preserves (in both directions) the atomic $\LL$-formulas, the standard (model theorists) answer to the so called {\em doxological challenge} is the following:

\begin{procedure}[Identifying an isomorphism type]\label{answer}
\begin{enumerate}[a)]
 a) Any particular structure is defined/described by its construction 
 \cite[p 9]{ButtonWalshbook}. E.g., the finite Von Neuman  
and   Zermelo ordinals are distinct but isomorphic structures (See Definition~\ref{omdef} for such a construction.). b) An isomorphism type of a structure $M$ is the collection (proper class) of all $N$ with  $N \approx M$.
\end{enumerate}
\end{procedure}

We argue here that the resources necessary to prove that a structure exists suffice to
pick out its isomorphism type. We do not find answer  b) of Procedure ~\ref{answer} in \cite{ButtonWalshbook}. They \cite[p. 151-2]{ButtonWalshbook} first distinguish a moderate objects modelist from an immoderate modelist on the criterion that the first denies a special faculty of `mathematical intuition'. Without the aid of such intuition, they suggest that a moderate modelist must reject `reference by acquaintance' and `must make do with reference by description', i.e. they must answer the doxological challenge `by laying down some formal theory'.  Taking  this second path they interpret the  formal theory as
a categorical theory and conclude that it cannot be first order by  L{\"o}wenheim-Skolem  considerations. No argument is given for the shift in meaning from {\em structures} to {\em models of a formal theory in a specific logic}, but in the context of page 152, it can be seen as an additional condition that must be met by a moderate-modelist. 

There seems to be no other option in \cite{ButtonWalshbook} aside from   invoking intuition or laying down a formal theory. However,  Procedure~\ref{answer}  provides such an option: a construction of the structure  in set theory (ZFC, though much weaker ones suffice, see \ref{omdef})  without making meta-theoretical investigations, not a theory in the same vocabulary as the structure whose models include the structure (and only those isomorphic to it). As we already mentioned, categoricity is a different notion altogether. It is a property, not of a structure or an isomorphism type, but of a theory $T$ in a specific logic: asserting there is only one isomorphism type of models of $T$. It is true that one could pick out an isomorphism type by finding a logic,  a theory,  and a structure $M\models T$ such that $T$ is categorical. But we are visiting Robin Hood's barn and the distracting trip around is 
more complex than the straightforward answer provided by Procedure~\ref{answer} in naive set theory.
However, the modelist has a more difficult task that we clarify with the discussion of frameworks in Section~\ref{carling}.

 One may say that Procedure~\ref{answer} is a mathematician's answer, one who overlooked a key word, the `{\em modelist's} doxological challenge', i.e. the challenge is for those who want to {\em explicate} mathematical structures by model theoretic tools. This recognition also clarifies the word `doxological', which refers to the fact that the explication is not intended as an {\em ontological} hypothesis (\cite[p. 38]{ButtonWalshbook}). This means that the challenge is not an epistemological question (\cite[pp. 145, 149]{ButtonWalshbook}), but one that requires an answer for how we form mathematical beliefs\footnote{\cite[p. 145]{ButtonWalshbook} emphasize that model theory is not of much help in answering epistemological questions concerning mathematical knowledge, but it is rather curious that someone would expect an answer from a modelist concerning the way we form mathematical beliefs. The answer of the modelist must certainly go beyond the resources of model theory.}. Can a mathematician be part of the modelists' community? As \cite[pp. 143, 164]{ButtonWalshbook} acknowledge, modelism is `an attitude towards model-theory' and `no one ever {\em called} themselves a `moderate modelist''. Thus, a practicing mathematician who has in his moments of reflection on his work this sort of philosophical attitude towards model-theory falls willy-nilly under this modelist label. 
 

One might object to Procedure~\ref{answer}:  `Set theory has raised its ugly head.' Answer: It already has\footnote{Shapiro acknowledges this,  \cite[p 181]{ButtonWalshbook}, in describing the `holistic' attitude to second-order logic.}. The formal definition of a structure requires a  weak set theory (neither the axioms of replacement nor choice).  Here is a canonical example. 

\begin{definition} \label{omdef} [The structure $N = (\omega,S)$]
Let $\tau$ contain a single unary function $S$. Let $S_N$ denote the function $\{\}$ mapping $x$  to $x \cup \{x\}$. 
\end{definition}
Halmos \cite[p. 51]{Halmos} uses induction to eliminate the possibility of finite cycles giving the existence of an $\omega$-sequence. \cite[\S 7]{Kunen} gives the standard account. He earlier lists  replacement as an axiom; it is not needed for existence of $\omega$;  that is provable in the weak Maclane set theory.  The axiom of choice is required to show 
Dedekind finite is equivalent to `embeddable in $\omega$'.
Maclane set theory has many equivalents; see \cite{Mathiasmst,Shulmansur}  or the summary in \cite[\S 3.2]{BaldwinMaddy}.





It should be clear by now that the original doxological challenge is a) distinct from any notion of categoricity and b) it is easily met by a (set-theoretic/logical (\cite{Reck})) modelist. In the next section, after we   
clarify the difficulties that confront the modelist in referring to a structure  by  specifying the different frameworks involved,  we show she can give a stronger answer to the
Doxological Challenge via categoricity.
While categoricity is   `out of reach for the first order theories we most care about' \cite[145]{ButtonWalshbook} (those  with finitary inference rules), moving to second order logic is an unnecessary step since by using the weaker $I-\omega$ and $G-\omega$-logics one regains the goal of characterizing canonical structures by categoricity.

\subsection{The Doxological Challenge and Carnap's Linguistic Frameworks }\label{carling}

Carnap’s (\cite[\S2]{Carext}) account of {\em linguistic frameworks} makes a fundamental distinction between two kinds of existence questions that one can formulate. To speak meaningfully about a new class of entities —numbers, properties, or propositions— one must first introduce a corresponding linguistic framework, that is, a system of expressions and rules governing the use of these expressions. Within such a framework, {\em internal questions} concern the existence or identity of particular entities as determined by the framework’s rules. These are legitimate theoretical questions, answerable by logical or empirical methods depending on the nature of the framework. By contrast, {\em external questions}—those asking whether the entire system of entities is `real' or `exists'—are, for Carnap, not genuine theoretical questions but pragmatic choices about adopting a mode of speech or a theoretical framework as such. They lack cognitive content and are therefore pseudo-questions.

Carnap illustrates this distinction with examples of 
frameworks for: things, natural numbers, propositions, thing properties, the integers and rationals, the real numbers, and the spatio-temporal coordinate system for physics. An internal statement such as `these two pieces of paper have at least one color in common' is an empirical claim within the framework of thing properties. By contrast, the metaphysical question whether {\em properties themselves} are real exemplifies an external question and, in Carnap’s view, it has no factual or cognitive significance. Similarly, within the physical framework, empirical questions about spatial or temporal relations are internal, while questions about the `reality' of space and time themselves are external pseudo-problems. In the framework for natural numbers one may likewise ask an internal question as `Is there a prime number greater than a hundred?' and this will be precisely answered by analysis based on the rules for the arithmetical expressions. However, if someone asks `Are there numbers?', as a question concerning the {\em ontological status of numbers as independent entities that have a reality of their own}, then this question simply lacks cognitive content for the natural numbers framework (\cite[\S2]{Carext}). Carnap's discussion of this last example illustrates very well the peculiarities of the external questions:

\begin{quote}Therefore nobody who meant the question ``Are there numbers?" in the internal sense would either assert or even seriously consider a negative answer. This makes it plausible to assume that those philosophers who treat the question of the existence of numbers as a serious philosophical problem and offer lengthy arguments on either side, do not have in mind the internal question. And indeed, if we were to ask them: ``Do you mean the question as to whether the framework of numbers, if we were to accept it, would be found to be empty or not?" they would probably reply: ``Not at all; we mean a question prior to the acceptance of the new framework." They might try to explain what they mean by saying that it is a question of the ontological status of numbers; the question whether or not numbers have a certain metaphysical characteristic called reality (but a kind of ideal reality, different from the material reality of the thing world) or subsistence or status of ``independent entities." Unfortunately, these philosophers have so far not given a formulation of their question in terms of the common scientific language. Therefore our judgment must be that they have not succeeded in giving to the external question and to the possible answers any cognitive content. Unless and until they supply a clear cognitive interpretation, we are justified in our suspicion that their question is a pseudo-question, that is, one disguised in the form of a theoretical question while in fact it is a non-theoretical; in the present case it is the practical problem whether or not to incorporate into the language the new linguistic forms which constitute the framework of numbers.\end{quote}

We consider the doxological challenge raised by Button and Walsh in relation to four {\em nested frameworks}, i.e. contexts that contain expressions for different categories of objects and rules that govern their use. Our basic objection to the Button-Walsh doxological challenge `How can we pick out particular
isomorphism types?' is exactly that given by Carnap for questions of the `reality of properties or numbers'. It is devoid of {\em cognitive content}, if frameworks are ignored  or  it is asked in an inappropriate framework. We reconstruct the study of mathematics as taking place in three different frameworks and address the Challenge in relation to them in the fourth framework, the one in which the philosophical reflection takes place.


 \cite[p. 152]{ButtonWalshbook} make a distinction between `reference by acquaintance' and `reference by description'. The first one is supposed to be obtained with the aid of mathematical intuition, but since the moderate modelist denies a special faculty of mathematical intuition, she aims to `pick out' reference by description. 

 We find the use of the word `reference' problematic in the context of mathematics and we consider that its use generates misunderstandings since people tend to use this expression as standing for a sort of `object'. For us, `mathematical objects' are concepts of different levels of complexity (depending on the framework to which they belong) that are constructed and grasped by the human mind. In each of the four frameworks below people who work in those frameworks are acquainted with and can describe the mathematical concepts from that framework.  

\begin{enumerate}
\item {\em Structure framework}. 
This is the  informal (natural language) framework in which children -- and likewise practicing mathematicians 9--  
consider intuitive or naive mathematical concepts, such as numerals, functions or operations that children work with. Children grasp these concepts in the process of learning (and so are acquainted with them) and can describe them by using the means of natural language. In particular, there is a naive notion of the natural numbers endowed  with functions ($+,\times$) and relations ($<$).
We refer to concepts at this level as {\em naive}.

\item {\em Theoretical framework of naive set theory}. This framework is employed in ordinary mathematical practice by mathematicians when they construct and compare various structures; in particular when they identify an isomorphism type.\footnote{The number theorist Barry Mazur \cite[p.222]{Mazur} describes very nicely the work of the practicing mathematician in this framework: `To define the mathematical objects we intend to study, we often -- perhaps always -- first make it understood, more often implicitly than explicitly, how we intend these objects to be presented to us, thereby delineating a kind of {\em super-object};
that is, a species of mathematical objects garnished with a repertoire of {\em modes of presentation}. Only once this is done do we try to erase the scaffolding of the presentation, to say when two of these super-objects --possibly presented to
us in wildly different ways -- are to be considered {\em equal}. In this oblique way, the
objects that we truly want enter the scene only defined as {\em equivalence classes of explicitly
presented objects}. That is, as specifically presented objects with the specific
presentation ignored, in the spirit of ``ham and eggs, but hold the ham." ' In the paper he suggests category theory rather set theory as the most suitable organization or scaffold (See \cite{BaldwinMaddy} for `ignoring').} Mathematicians define in naive set theory
\footnote{ We use the term `naive set theory' in the sense of \cite{Halmos}, or the mathematician who is 
unfamiliar with any details of axiomatic set theory.} both individual structures (an $\omega$-sequence) and classes of structures in naive set theory (the class of semigroups, structures $(A,+)$ such for all $a,b,c\in A$, $(a+b) +c = a+ (b+c)$). Note there is no formal language present. The first corresponds to what \cite{Sha97} calls a non-algebraic theory (i.e. univocal theory) and the second to an algebraic theory. We avoid `theory' in this context, reserving it for formal theories. We suggest that mathematicians have {\em acquaintance by definition/construction} to  structures defined in this manner. Following \cite[p. 9]{ButtonWalshbook}  we refer to the concepts used in this framework as {\em informal}.  Some mathematicians might prefer `formalizable' as they are consciously
working loosely in ZFC.

\item {\em Metatheoretical framework of formalized theories in various logics}. This framework contains formal mathematical concepts defined by the mathematical logician. She is acquainted with them by defining a formal axiomatization of the informal concepts from Framework 2 and can describe some of those concepts by a formal logic, a vocabulary and a categorical set of axioms. 
More specifically, those who work in this framework are directly acquainted with the {\em formal concept} and can refer by a {\em formal description}
\footnote{As \cite[p 152]{ButtonWalshbook}  require of a moderate or conceptual modelist.} to the concept from Framework 2 that is categorically formalized.

 \item{\em Philosophical framework.}
This is the framework in which one can
reflect on any philosophical features arising in the first three. Sections \S \ref{Dochall} and \S \ref{AgCA} of this paper take place in this philosophical framework. For instance, in this framework we can also state the following principle about the previous three frameworks:
 
 {\em Meta-framework principle}: If a question is properly asked and answered in framework $n$, then it remains so in higher frameworks.

\end{enumerate}

 The Doxological Challenge gets a clear answer once we set up a precise and adequate framework, and is meaningless otherwise. Consider Frameworks 1 and 2. As children learn the informal nested framework of things, properties, and natural numbers and understand that a structure contains functions and relations, then they can talk about that particular structure, i.e. they use the `structure framework'. In this framework children refer to natural numbers by simply using the language and they can answer simple (internal) questions. It is meaningless to ask them whether these numbers `really' exist. The mathematician can answer the same internal questions, but he can also describe what his arithmetical language refers to once he has adopted the theoretical framework of naive set theory.  In the latter framework of set theory the mathematician can simply refer to a particular structure due to the fact that he actually constructs it. This is the reason that the Doxological Challenge gets a clear and precise answer by Procedure~\ref{answer} in Framework 2.

However the modelist has a more difficult task, which requires Framework 3 (and Framework 4 for the discussion).  If someone wants to formally represent a certain structure, she has to know   (i.e. have a concept of) what she wants to formally represent.  The model theory approach provides a common ground for all signatures and many (all?) logics and even for the formalism-free approaches to model theory. The categoricity approach requires a different approach and a possibly different logic to describe the isomorphism type for each structure.

The move to Framework 3 may  both sharpen and extend the mathematical conception of a rougher notion. 
For example, Euclid, with axioms that are informal in our sense,
 imagines a line with some irrational numbers described by specific algorithms; Dedekind's 2nd order logic approach imagines a line with irrational numbers described by arbitrary `cuts'. Both are working in rudimentary versions of Framework 3. While Dedekind fits well\footnote{Certainly, if we adopt Quine's view of second order logic.} in Framework 2, by stressing  definitions and axioms Euclid presages Framework 3. Framework 3 has allowed model theorists to make deep contributions to analysis, number theory, differential equation, combinatorics etc. Our approach in proving categoricity with the $I-\omega$ and $G-\omega$-logics was conducted in this framework.

 The categoricity, decidability, completeness, et al. questions are internal questions in Framework 3, where formalized theories in various logics are analyzed, but these questions are meaningless (external) questions outside this metatheoretical framework (for instance in Frameworks 1 and 2). In particular, the question about uniqueness of {\em models} up to isomorphism and categoricity of formal theories are internal questions in metatheoretical frameworks of formalized theories in a certain logic, but external to Framework 2. The questions about the identification and uniqueness of {\em isomorphism types} are internal questions in both Frameworks 2 and 3. The difference is that in Framework 2 the mathematician identifies the isomorphism type by constructing/defining it, while in Framework 3 the mathematical logician identifies it by choosing an adequate logic to obtain a categorical theory that it satisfies.\footnote{We note in passing that when someone provides derivations in a formal theory in a given logic, he does not need the whole resources of Framework 3. This is a more narrow framework that comprises the axioms of the theory and the deductive apparatus of its underlying logic.} This step from Framework 2 to 3 introduces formalized
 mathematics and is the answer by the moderate modelist to the doxological challenge. It provides {\em reference by (formal) description} in the sense of \cite{ButtonWalshbook}. In Section~\ref{AgCA} we show that the construction of our inferential $\omega$-logics does not invoke the concepts aimed to be categorically formalized.





 
 In discussing the doxoloxical challenge, our emphasis is on the distinction between the first two frameworks: the Challenge is external to the structure framework and easily answered in the framework of naive set theory. Categoricity of theories is simply irrelevant in these two frameworks since it presupposes a theory formalized in a certain system of logic.\footnote{ In Framework 3, one can consider the question of categoricity of theories in a particular logic. The question is external to the logic itself, where one works in the object language of that logic; it is addressed and answered in the metatheoretical framework of set theory. Likewise, in this framework one can compare the resources needed and the results provable or at least  stateable in the logic. The assertion that the continuum hypothesis is true, for example, can be formulated for $PA^2$ but not for weaker arithmetical systems. But questions motivating such comparisons arise in Framework 4.}

\subsection{Cognitive Modelism}\label{cogmodelism}

 Our own philosophical view is closer to what  \cite{ButtonWalshbook} call {\em concepts modelism} since we see classical mathematics as a complex process of constructing and developing a certain kind of concept. The claim in this sentence is  not intended to be an ontological one: when we say `structure' in this paper we  have in mind something that has a conceptual nature, not an object that inhabits Plato's realm. We  refer to  concepts defined and grasped by human mind, whose linguistic expression\footnote{For practical purposes, we shall not defend here a systematic philosophical conception on concepts. For a general discussion on concepts see \cite{MargolisLaurence} and their \S4 for the relation between language and concepts.} is a constitutive component. We do not make the further step in postulating a realm of objects that correspond to these concepts and we see this as a problematic step that blunt objects modelists make in accordance with their {\em metaphysical credo}.  That certain mathematical structures defined by mathematicians are suggested by configurations found in the physical world or that some mathematical concepts are sometimes found to be partly instantiated in the physical world does not undermine the conceptual nature of the mathematical structures in our view. Mathematical concepts are likewise important for objects modelists since, for instance, the assertion that there actually are $\omega$-sequences presupposes that we possess the concept of $\omega$-sequence.
 
 \cite{ButtonWalshbook} use terms such as `natural number structure' to designate a particular isomorphism type. However, this is 
meaningful just if there is a concept of natural number structure. \cite[p. 150]{ButtonWalshbook} discuss how humans acquire, express, share such concepts as redness and write `In the case of redness, though, we can point to both red things and also to organs which are good at detecting red things. In the case of mathematical concepts, such as {\em $\omega$-sequence}, it is doubtful that there is anything similar to point to.' 
We adopt the name, {\em cognitive modelism} given by
\cite{DeBenRos}, for their description of `what to point to'  while refining the notion using the frameworks established in Section~\ref{carling}.

The quote above from \cite{ButtonWalshbook}  oversimplifies the development of the concept `red'.
A long time ago people  identified certain properties of objects as colors
and others as shapes.  Much later\footnote{Likely in the last 2000 years; an internet glance at `red  in Indo-European' showed a commonality - the first letter was usually `r'.}, a particular group of people decided a certain family of colors
that were not too different should be called `red'.  Of course there are disagreements over
whether a particular scarf is red or scarlet. Eventually, an international agency decides on
a certain family of frequencies.  Mathematical concepts evolve in a similar way.

\cite{DeBenRos} give an account of how children acquire the concept of an $\omega$-sequence. Of course, the specific organ(s) must include the brain. But rather than merely asking, `What color is the ball?', the study of concept acquisition requires much more complicated analysis. \cite{DeBenRos} partially base their account on the work
of such scientists as \cite{Carey} analysis of children's development
of number concepts\footnote{In addition to the bibliography in \cite{DeBenRos} consult such authors as \cite{Devlingene,LakNun}.}. The book,  \cite{Devlingene}, argues extensively 
with wide ranging references to cognitive science: `the features of the brain that allow us to do mathematics are the very same features that enable us to use language
-- to speak to others and understand what they say'.

While the early steps of this development are known through observations of
babies, the process continues through school. But even more, it continues through history. Although Euclid was well-aware of ratio and proportion, 
he admitted neither the unit nor rationals as numbers \cite[\S 2.1,\S 3.1]{Mueller}. Today, we maintain
the second exclusion in our notion of `natural' number, but not the first.
In a non-empirical direction, \cite{Whiteheadpr} 
analyzes how the mind moves from sensory perception to
abstract concepts.

There are other ways of grasping a structure than having a picture in mind of its quantifier free diagram, which works for $\omega$-sequence.  Grasping the structure of natural number arithmetic by
`visualizing the diagram' is problematic.  One only knows the addition and multiplication by rules, there are sentences of arbitrarily high quantifier complexity, and easily stated but difficult to picture, let alone prove, such as Fermat and twin primes.

 The practice of mathematics illustrates very well the fact that mathematical concepts rarely have a rigorously defined content from the very start. The process of {\em concept formation} is dynamic and a categorical formal concept is an advanced stage in the process of making a mathematical concept precise. It is thus natural to invoke some pre-formalized/informal mathematical concepts in laying down a (categorical) formal theory in a certain logic. Even in formulating the semantics of first-order logic we use notions as `interpretation', which relies on a previous grasp of the notion of `function'.


On our view, there is a concept of natural number structure\footnote{\label{fnamb} Actually, there are several. In a footnote to \cite[\S 2]{ButtonWalshart}, but not in the book, Button and Walsh are clear to distinguish $\omega$-sequence, vocabulary: $\langle 0,S\rangle$, from arithmetic, $\langle 0,+ ,\times\rangle$. That is crucial, as it is only with some sort of induction axiom (schema) that one can define  the arithmetic operations on an $\omega$-sequence. Both concepts appear in Framework 1 and are grasped by most adults, though perhaps not differentiated.}. That is partly grasped by those who work in Framework 1 and fully grasped by those who work in Framework 2. The goal of categoricity is achieved by choosing a vocabulary, a theory and a logic in that vocabulary, specified by its rules of inference, and axioms  such that all models of those axioms are isomorphic. We have provided two such logics and appropriate rules of inference and axioms. Certainly, those who work in Framework 3 also grasp a concept of natural number structure and simply aim to formally capture it by means of a categorical theory.

\cite{DeBenRos} argue that the standard natural numbers are to be found in the concepts that cognitive agents acquire in a learning processes. In particular, they reconstruct the cognitive process as modeled by a function between numerals and suitable collections corresponding to them that is generalized by induction to the entire number sequence. This rational reconstruction of the learning process concludes that cognitive agents do grasp the concept of an omega sequence. We are sympathetic with this approach and we take the evidence provided by Cognitive Science as a very reasonable justification for the idea that in the learning process children become acquainted with {\em the intuitive} notion of an $\omega$-sequence. This idea in turn justifies the intelligibility of the use of the $\omega$-rule in our approach since grasping its premises will be no mystery even for children, who will quickly understand its premises as asserting that each object in the sequence they learned has a certain property.\footnote{\cite[Ch. 10]{Warren} also argued that the categoricity of $PA$ can be obtained by using the classical generalized $\omega$-rule as an {\em open-ended} rule schema (i.e. a rule that remains sound in extensions of the initial language/vocabulary). A non-standard model is blocked by considering an instance of the rule with the one-place predicate `$SM(x)$' that means `x is a standard number of model $M$'. The conclusion of this instance (i.e. $(\forall x)(N(x) \rightarrow SM(x))$) will be false in a non-standard model $M$ and, thus, the model will not be admissible. Mathematically, we don't find intelligible the use of this predicate in an instance of the rule as applied in arithmetic. It may be attractive as a philosophical argument (in Framework 4), but not as means for categorically formalizing $PA$ in Framework 3 (see also \cite[\S 7]{Brincus} for discussion). In particular, since categoricity is a property of theories formalized in a certain logic, it must be achieved by mathematical means in Framework 3, i.e. {\em there is no philosophical substitute for mathematics}. On the philosophical side, Warren's argument for the possibility of following the $\omega$-rule is interesting and we are sympathetic to it (see \cite{Topey} for a recent critical discussion). We agree with \cite{Church} that `it is clearly not possible for the users of a language systematically to follow a non-effective rule in practice', and that these  rules `have a place, and perhaps an important one, in theoretical syntax.', i.e. in Framework 3.} This may also explain Tarski's \cite[p.411]{Tarskilogcon} natural attitude of taking the classical $\omega$-rule as {\em intuitively} valid.

Unlike the rules of inference, the concepts required for axiomatization
depend on the topic and lay out in the vocabulary the fundamental notions being axiomatized. This explanation is lost in second order logic as quantifying over concepts allows the formal definition of virtually anything conceivable.

As elaborated in \cite{Baldwinphilbook}, the vocabulary is chosen based on the informal concept held by the formalizer. 
Two informal concepts are often conflated as `natural numbers', but properly distinguished in \cite[\S 2]{ButtonWalshart}: $\omega$-sequence and arithmetic. The conflation occurs because of the immense power of second order definition \cite[\S 5.B]{ButtonWalshbook} or even first order induction schema. We consider first order axioms for the examples in ({Theorem~\ref{wrulebworks} 
with specified vocabulary for $\omega$-sequence as  an $(S,0)$ structure and for arithmetic $(S,0,<, +,\times)$.} With our  (non-circular) $I-\omega$-rules we establish categoricity of the appropriate theories. 

However, the simplest first order axiomatization of an $\omega$-sequence is to require that $S$ is $1$-$1$ and onto every element
except $0$
(\cite[Definition 1.9: Q1-Q3]{ButtonWalshbook}) and require that there are no finite cycles. This gives an first order axiomatizable $\aleph_1$ (and so 
uncountably)-categorical structure whose prime model is $(\omega,S)$. The structure is even strongly minimal (every definable subset is finite or cofinite). Such definable sets
are the building blocks of theories categorical in uncountable cardinalities. 
This is the `trivial' case of the Zilber trichotomy (below). 
That is, far from being a wild theory such as arithmetic, the unique structure given by the $I-\omega$-rule and these axioms is an iconically simple decidable theory. On the other hand, adding the induction schema so as to define  addition and multiplication   \cite[Definition 1.9]{ButtonWalshbook} gives full first order Peano Arithmetic.

\begin{remark}[Three examples of mathematical concept formation]\label{ZilHru}{\rm

The concept of countability\footnote {That is, able to be counted. So mathematically, it means finite or of cardinality $\aleph_0$.} has a long history.  Aristotle distinguished between `potential' and `completed infinity' and rejected the latter. The notion of various `orders of infinity' appears at least by  the 12th century writing of Grosseteste \cite[p 134]{Grosse}.
But the modern fully worked out concept of comparing infinities arrives only with \cite{Cantor}.\footnote{\cite{Kajercline} has recently argued that the meanings of the numerals are essentially tied to counting. We agree with this idea as an explanation of what happens in Framework 1, when children learn numerals, but we see no conflict between this idea and the approaches conducted in Frameworks 2 and 3.}

 \cite{Giaquinto}, argues that a structure/isomorphism type can be known `more directly than as the structure of all models' of a certain categorical theory. He lays out how a cognitive agent can visualize finite structures, abstract to `templates' and then organize sequences of templates to
obtain knowledge of a structure, as that of  natural numbers. This knowledge has an experiential constituent and it is thus different from `knowledge by a  description of the form {\em the structure of such-\&-such axioms'} \cite[p 57]{Giaquinto}, which requires knowing that the formally axiomatized theory is categorical.

Here is a modern example of developing the concept of a specific structure.
In proving in the 1980's that there is no finitely axiomatizable first order theory that is categorical in both $\aleph_0$ and $\aleph_1$,
Zilber introduced a classification of strongly minimal theories.
Strongly minimal theories are the simplest type of theory
that is categorical in all uncountable powers.
A strongly minimal theory admits a combinatorial geometry induced
by the notion of algebraic closure ($a \in \acl(B)$ if for some
$\phi(x,\bbar)$ with $\bbar \in B$, $\phi(a,\bbar)$ and $\phi(x,\bbar)$ has only finitely many solutions).
Roughly speaking, the geometry is discrete  
if $a\in \acl(B)$ implies  $a\in \acl(b)$ for a single element of $B$
and locally
modular if the closure relation behaves like a vector space. Otherwise the closure is non-locally modular. 
Zilber conjectured that every non-locally modular geometry
of a strongly minimal set behaves like a field.

In the early 1990's Hrushovski refuted this conjecture by defining a class of finite structures
with an axiomatically defined strong substructure relation\footnote{Generalizing, the Fra{\" \i}ss\'{e} limit (a 1950's construction) that  moves from
a uniformly locally finite class of finite structures to a countable $\aleph_0$-categorical
structure via substructure.}, $\leq$.
The  direct limit under $\leq$ of the members of the class is a
is a countable structure. It is
strongly minimal, not $\aleph_0$-categorical, not locally modular and not field-like. Model theorists grasp this structure.
Hrushovski used a similar technique to refute the conjecture
that every stable $\aleph_0$ categorical theory is superstable.
These methods have been used in hundreds of papers since. These papers adapt the construction to
find theories in many classes of the classification hierarchy.   There is no
need for being able to visualize the structure, although some
extremely geometric examples were constructed with a partial visualization in mind. The moral of the story is that an explicit  description of the structure is not essential.

It is by no means clear that a Hrushovki theory has an algebraically prime model that would make it categorical
in $I-\omega$-logic. However, the Scott sentence of the saturated countable model
in a vocabulary with constants naming a basis for the geometry,
is non-generative and so it is categorical 
in $G-\omega$-logic.\footnote{The second order theory
of an infinite set of constants is not categorical.  The first paragraph
of \cite{Kesk} begins with a short argument to this effect by Shelah.}
}
\end{remark}

\section{Against the Circularity Argument}\label{AgCA} 
\stepcounter{subsection}

In this section we address what we call the `circularity argument' of
Button and  Walsh.  We work here in Framework 4, since we are discussing philosophical issues about the Framework 3  analysis of formal theories.  
They  \cite[p.164]{ButtonWalshbook} argue that 
an attempt by a moderate modelist  to appeal to a categoricity theorem in an extension of first order logic
to explain  how we grasp a mathematical concept `by formal description' is hopeless. As,  one must   semantically set forth a logic  that   `will invoke precisely the kind of mathematical concepts that were at issue in the first place, and which they were hoping to secure by appeal to a categoricity theorem'.  We rebut this argument by showing categoricity in an infinitary inferential logic.


This circularity argument is also described by \cite[p.151]{ButtonWalshbook} as a reiteration of the Doxological Challenge at the level of selecting a logic because the logics stronger than first order logic `must be characterized using mathematical concepts'.
 They don't specify what concepts are involved. Clearly finite, countable, and the inductive definition of formulas and satisfaction are not included as they are used  to formulate the syntax of first order logic.

We observe in this section that our inferential approach to categoricity
requires only the notions of finite and countably infinite and
so does not presuppose the notion of an $\omega$-sequence or arithmetical notions. For this we consider the precise definition of the notions used in formulating our inferential rules.

\begin{definition}\label{carddef} We say two sets have the same {\em cardinality} if they can be put in $1$-$1$ correspondence.
 Following Cantor \cite[p 86]{Cantor}  and Dedekind, 
a set $D$ is {\em (Dedekind finite)} iff every one-one map of $D$ into itself is onto. We take this as our notion of finite\footnote{It is well-known that the equivalence of 
Dedekind finite with the standard definition in e.g. \cite{Kunen}, $|A| < \omega$,
where $|A|$ is defined as the least ordinality of a set in $1$-$1$ correspondence with $A$
requires the axiom of choice. For the moment, we disregard this search for unique reference that is not necessary for our purposes.}\label{cardstr}.
Further, we say a set is {\em infinite} if it is not finite.  Finally a set $X$ is countable
if it is infinite and every proper subset is either finite or has the same cardinality as $X$.
\end{definition}

Definition~\ref{carddef} is in
in Framework 2. Thus the notions of (equi)-cardinality and countable are available to construct logics in Framework 3.

\begin{remark}{Our rules are based on cardinality of at most countable sets,
 not arithmetic}.\label{cardrules}
{\rm Recall the sketch of the history of the concept `countable' in Remark~\ref{ZilHru}. Note that our $\omega$-rules, (\ref{ourinfrule}, \ref{ourinfrule2}), all have the form: if a formula holds for all constants in a countable set $C$, then the universal closure of the formula also holds. Thus, our rules of inference\footnote{ Definition~\ref{Morcode} uses a numbering of the constants that smells of arithmetic (although only a pairing function). But this is used in the proof of Theorem~\ref{omegaruleatomic2}, not in the 
formulation of the rule.}
depend only on the notion of finite and countable cardinality. Neither the rules of inference, nor the truth conditions for infinite conjunction
or disjunction have any dependence on the crucial concepts
of well-order and induction aimed to be secured by the categoricity theorem for $PA^\omega$. 
Moreover, this approach extends to $L_{\omega_1,\omega}$ via
Theorem~\ref{omegaruleatomic2}.}
\end{remark}

We now contrast our logics with other weak logics discussed in 
quantifier (\cite[p 163]{ButtonWalshbook}). 
\cite[\S 7.10]{ButtonWalshbook} do not analyze
$\omega$-logic, because it has an infinitary rule of inference.  But they argue that each of seven logics weaker than second order are `just more theory'\footnote{ \cite{Speitel2}, especially \S 2.3 on `just more theory' provides an argument very similar to ours.
By distinguishing the construction of the structure in Framework 2
from the categoricity in Framework 3, our argument is more general. Our argument for the concepts of cardinality and infinity being  formalizable in 
Framework 2 justifies the semantic characterization
of Speitel's logic $L(Q_0)$ (as well as the H\"{a}rtig equicardinality quantifier) and
authorizes the logics, not just the categoricity of arithmetic.}. They consider that all these weaker logics must be semantically specified and so are dismissed for invoking in their semantics precisely the concepts aimed to be `secured' by a categoricity theorem. In contrast to the previous seven, our logics are inferentially specified. 
 Moreover, our logics are no more semantically defined than first order; perhaps less, our $\omega$-rules yield categoricity of the quantifiers and the $\Rbar$-semantics are thus justified syntactically.

\cite[pp. 163-4]{ButtonWalshbook} argue  against $L_{\omega_1,\omega}$  that `grasping an infinite disjunction is just like grasping the natural number sequence'. We do not accept  this equivalence. Note that an infinite disjunction is derivable from any of its disjuncts and, on the semantic side, the truth value of an infinite disjunction depends on the set of formulas, not their order. In Theorem~\ref{omegaruleatomic2} we reduced 
$L_{\omega_1,\omega}$ to $G-\omega$-logic. As we argued earlier, each of our $I-\omega$ and $G-\omega$-rules depend only on a notion of a countably infinite infinite set of constants, not any relations that hold among those constants. And, on the notion that there is a point that is not in a given set.

 Koellner  (\cite[p.19]{Koellner}) argued that either we view the $\omega$-logic model-theoretically or proof-theoretically ``it is hard to maintain that results in $\omega$-logic depend on syntactic form alone”, since in either case this logic is strongly connected with the standard model of arithmetic. This argument simply does not apply to our inferential logics; Kollner works with the historical notion of $\omega$-logic that concerns only theories that interpret arithmetic and uses the notion of $\omega$-model in 
Remark ~\ref{omega-logics and rules}.2.a.i. Our inferential
$\omega$-rules apply to any first order theory -- not just to those interpreting first order arithmetic, and the corresponding notions of model for these logics are those from Remark~\ref{omega-logics and rules}.2.a.ii  and  Remark~\ref{omega-logics and rules}.2.b. Although there is a strong historical connection between the $\omega$-rule and the standard model of arithmetic,  our $\omega$-rules use only the notions in Definition~\ref{carddef}. 
 
 In building up our inferential $\omega$-logics we have added to the syntactical instruments  constructing  first-order logic by the restriction to  a countable number of individual constants -an idea that is common in the discussion of logical inferentialism since \cite{CarnapFor}. Philosophically speaking, we say that we have assumed `nameability' and `countability'.
 
  We noted above (Remarks~\ref{ZilHru} and ~\ref{carddef}) that the notion of `countable' does not rely on arithmetic. Nameability is a constitutive feature of a natural language, in the sense that ordinary speakers use many expressions with the intention of naming something. In addition, it is usual linguistic practice to introduce a new name in the language when  a new object enters into the universe of discourse. Although nameability may be seen as a model-theoretic or referential assumption, we take it to be compatible with an inferentialist perspective, since introducing a name for something does not immediately provide the name with a meaning. The meaning of that name is determined by
  by the rules that govern its use in linguistic and inferential practices. Finally, nameability is certainly not an arithmetical concept.

Thus, the formalization of first-order theories in the $I-\omega$-logic and $G-\omega$-logic is, indeed, `more theory', but it is not `just more theory', since the categoricity problem gets a clear solution by theoretical tools that do not presuppose notions they aim to justify.

Does the choice of logic give more information than is needed or encode irrelevant questions?
Such immodesty\footnote{See \cite[\S 11.2]{Baldwinphilbook} for `immodesty of Dedekind's axiomatization of geometry.} is a sin
of second order logic; the Turing degree of 2nd order arithmetic is `fundamentally impossible to describe in a concrete way'
\footnote{See answer to \url{https://math.stackexchange.com/questions/3644483/what-is-the-turing-degree-of-the-set-of-true-formula-of-second-order-arithmetic}.}.
Many less expressive logics suffice to prove categoricity of arithmetic. While,  second order offers  the  categoricity of the reals as a bonus. And the continuum hypothesis has been decided. But, which way will depend on our choice of the ambient set theory as a metatheoretical framework.

We compare our approach with Dedekind's second order axiomatization of arithmetic. 

\begin{enumerate}
    \item We obtain categoricity in a first order logical framework, i.e. with quantification over individual variables.  
    \item 
    We univocally read off the
    truth-conditions of the quantifiers from our rules of inference. Second order logic (SOL) and most of the weaker logics use semantic notions. SOL needs the standard semantics for categoricity, while its rules of inference are complete for the Henkin semantics. 
  
    \item On the semantical side our quantifiers range over points as in first order logic, while second order arithmetic can define  arbitrarily complex relations \cite[\S 5.B]{ButtonWalshbook}. 
    So, ontologically, we are as conservative as possible.
    
\end{enumerate}

While the classical $\omega$-rule provides the negation completeness of $PA$ (\cite{CarnapLog},\cite{Shoenfieldomega}) but not its categoricity, our approach shows that with an inferential version of the $\omega$-rule, interpreted in terms of R-semantics, we obtain the negation completeness of $PA$ and its categoricity. In SOL we obtain only the second meta-theoretical property. With these criteria, a costs-benefits analysis  favors our approach. 



\section{Conclusion.}

 We argued in this paper that the mathematical results obtained by means of the inferential $\omega$-logics answer 
 the doxological challenge and the circularity argument raised by Button and Walsh \cite{ButtonWalshbook}. In particular, we argued that the doxological challenge is generated by conflating the structure and naive set theoretic frameworks described in Section~\ref{carling}. Likewise, we argued that the circularity argument does not apply to the proofs of the categoricity theorems in our $\omega$-logics since our $\omega$-rules do not presuppose arithmetical concepts that were in need of determinacy by a categoricity result.

\bibliography{cbjb.bib}

@incollection{BaldwinMaddy,
    author = {John T. Baldwin},
    title = {Exploring the generous arena},
    booktitle = {The Philosophy of Penelope Maddy},
    editors = {J. Kennedy and S. Arbeiter},
    publisher = {Springer, Outstanding Contributions to Logic},
    year = {2024},
pages = {143-164}
     }

@book{Baldwinphilbook,
author = " Baldwin, John T. ",
year={ 2018},
publisher ={Cambridge University Press},
title = " Model Theory and the Philosophy of Mathematical Practice:
 Formalization without Foundationalism",
PAGES = {352}
}

@article{BaldwinBrincusI,
  author  = {Baldwin, John T. and Br{\^i}ncu{\c{s}}, Constantin C.},
  title   = {Categoricity by inferential $\omega$-logic and {$L_{\omega_1,\omega}$}},
  year    = {2026},
  note    = {\url{https://arxiv.org/pdf/2602.02854}}
}

@article{Brincus,
title ="Categorical Quantification",
author ="Br{\^i}ncu{\c{s}}, Constantin C.",
journal ="The Bulletin of Symbolic Logic",
volume = {30},
pages ={227-252},
year ={2024},}

@article{ButtonWalshart,
title ="Structure and Categoricity: Determinacy of Reference and Truth Value in the Philosophy of Mathematics",
author ="Tim Button, Sean Walsh",
journal ="Philosophia Mathematica",
volume = {24:3},
pages ={283-307},
year ={2016},}

@book{ButtonWalshbook,
title ="Philosophy and Model Theory",
year = {2018},
author ="Button, Tim and Walsh, Sean ",
publisher ={Oxford University Press},
pages = {xvi+517}
}

@book{Cantor,
 author ="Cantor, Georg",
 year = "1915",
title = "Contributions to the founding of the theory of transfinite numbers", 
pages = "211", 
publisher = "Dover", 
 address = "New York",
 note ={edited by Philip E.B. Jourdain.  German publication 1895, 1897, this edition 1915}}

@book{Carey,
author ="Carey, Susan",
year = "2009",
title = "The origin of concepts",
publisher = "Oxford University Press",
pages = {598},
address = "Oxford"}

@article{Carext,
title ="EMPIRICISM, SEMANTICS, AND
ONTOLOGY",
author ={Carnap, R.},
journal ="Revue Internationale de Philosophie ",
volume ={4},
year ={1950},
pages={20-40},
note = {Reprinted in
the Supplement to {\em Meaning and Necessity: A Study in Semantics
and Modal Logic}, enlarged edition (University of Chicago Press,
1956)}
}

@book{CarnapLog,
title = "Logical Syntax of Language",
author = {Carnap, Rudolf},
publisher = "Trench, Trubner and Co Ltdr",
address ={London},
year =  "1934/37"}

@book{CarnapFor,
title ="Formalization of Logic",
author ={Carnap, Rudolf},
publisher ="Harvard University Press",
address ={Cambridge},
year = "1943"}

@book{ChangKeisler,
author ="Chang,C.C. and Keisler, H.J.",
year = "1990",
title = "Model theory",
publisher = "North-Holland",
note ="3rd impression 1992",
pages ="XII+550"}

@article{Church,
    author = {Church, Alonzo},
    title = {Review of {C}arnap 1943},
    journal = {The Philosophical Review 53:5},
    year = {1944},
    pages ={493-8}
}

@article{DeBenRos,
author = "De Benedetto, Matteo and  Rossi, Lorenzo", 
title ="Cognitive {M}odelism",
journal ="Philosophia Mathematica",
year = {2026},
note = 
"forthcoming" }

@book{Devlingene, title = "The Math Gene: How Mathematical Thinking Evolved 
and Why Numbers Are Like Gossip " ,
 year = " 2000" ,                              
 publisher = " Basic Book"                 ,
 author = "  Devlin, Keith" , 
pages = "xi + 425"}

@incollection{Ebbing,
author ={Ebbinghaus, H.D.},
title ="
Extended Logics: The General Framework",
booktitle ="Model-Theoretic Logics", 
editor = {J.~Barwise and S.~Feferman},
pages ={25-76},
publisher = {Springer-Verlag},
year ={1985}}

@incollection{Giaquinto,
author = " Giaquinto,M."    ,
year = "2008"           ,
title = "Cognition of structures ",
booktitle={The Philosophy of Mathematical Practice},
PUBLISHER ="Oxford University Press",
pages ={198-256},
editor ="P. Mancosu"
}

@book{Grosse,
author = "Anne Freemantle",
year = "1954",
title = "The Age of Belief",
publisher = {New American Library of World Literature},
series = {The Mentor Philosophers},
note = {Grosseteste quoted page 134}
}

@article{Kajercline,
title ="Don’t count on structure",
author ="Kajercline, Hayden",
journal ="Philosophical Studies",
year ={2025},
note ={https://doi.org/10.1007/s11098-025-02432-7}
}

@article{Kesk, 
title ="Characterizing all models in infinite cardinalities", 
 author ="Keskinen,Lauri ",
journal = "Annals of Pure and Applied Logic", 
volume = {164}, 
 pages={230-250}, 
 year ={2013} }

@article{Henkinomt,
author = " Leon Henkin",
title = "A generalization of the notion of $\omega$-consistent",
JOURNAL = {Journal of Symbolic Logic},
volume = {19},
pages = "183-196",
year = "1954"}

@unpublished{Koellner,
title = "Carnap on the Foundations of Logic and Mathematics",
author = "Koellner, Peter ",
year = "2009",
note = {\url{https://bpb-us-e1.wpmucdn.com/websites.harvard.edu/dist/f/94/files/2022/07/CFLM.pdf}}
}

@book{Halmos,
title = "Naive Set Theory", 
year = " 1960"                               , 
publisher = {Van Nostrand},    
author = "Paul  Halmos",
 pages ={vii +104 } }

@book{LakNun,
  title     = {Where Mathematics Comes From},
  author    = {Lakoff, George and N\'u\~nez, Rafael},
  year      = {2000},
  publisher = {Basic Books}
}

@book{Kunen,
title = "Set {T}heory, {A}n {I}ntroduction to {I}ndependence {P}roofs",
year = " 1980"                               ,
publisher = "North Holland",
author = "K. Kunen"}

@InCollection{MargolisLaurence,
	author       =	{Margolis, Eric and Laurence, Stephen},
	title        =	{{Concepts}},
	booktitle    =	{The {Stanford} Encyclopedia of Philosophy},
	editor       =	{Edward N. Zalta and Uri Nodelman},
	howpublished =	{\url{https://plato.stanford.edu/archives/fall2023/entries/concepts/}},
	year         =	{2023},
	edition      =	{{F}all 2023},
	publisher    =	{Metaphysics Research Lab, Stanford University}
}

@book{Markerbook,
author ="Marker,D.", 
title ="Model Theory: An {I}ntroduction", 
publisher ="Springer-Verlag", 
year ="2002" }

@article{Mathiasmst,
    author = {A.R.D. Mathias},
    title = {The strength of MacLane set theory},
    journal = {Annals of Pure and Applied Logic},
    year = {2001},
    pages = {100:107-234}
}

@incollection{Mazur,
    author = {Barry Mazur},
    title = {When is One Thing Equal to Some Other Thing?},
    booktitle = {Proof and Other Dilemmas: Mathematics and Philosophy},
    editors = {Gold B, Simons RA},
    publisher = {Spectrum. Mathematical Association of America},
    year = {2009},
    pages = {221-242}
}

@incollection{MorleyLeeds,
author = "Morley, M.",
title = "Partitions and Models",
booktitle = "Proceedings of the Summer School in Logic Leeds, 1967",
year =" 1968",
editor ={M.H. L\"{o}b},
note= {paper written by Vivienne Morley},
publisher = {Springer-Verlag},
addresss ={Berlin-Heidelberg-New York},
pages ="109-158"}

@book{Mueller,
author = "Ian Mueller " ,
title= "{Philosophy of Mathematics and Deductive Structure in Euclid's Elements}"  ,
series ={Dover Books in Mathematics},
year = " 2006",
PUBLISHER = " Dover Publications, Mineola N.Y." ,
note={First published by MIT press in 1981},
pages ={xiv + 378}
 }

@article{Orey,
author = " S. Orey"   ,
year = "1956"                          ,
title = "On $\omega$-consistency and related properties",
journal = "Journal of Symbolic Logic",
volume = {21 },
pages =" 246-252    "
}

@article{Reck,
    author = {H. Reck},
    title = {Dedekind’s structuralism: An interpretation and partial defense},
    journal = {Synthese},
    year = {2003}
}

@article{Shulmansur,
    author = {Michael Shulman},
    title = {Comparing material and structural set theories},
    journal = {Annals of
Pure and Applied Logic, 170:4},
    year = {2019},
    pages = {170:465-504},   
}

@incollection{Shapiro,
booktitle = "Handbook of Philosophical Logic",
author = {Shapiro, Stewart },
title = {Systems between First-Order and Second-Order Logics},
volume = {1},
year = " 2001"                               ,
editor = "D.M. Gabbay, F. Guenthner (eds) ",
publisher = "Springer",
address ={Dordrecht},
pages = "131-187",
}

@book{Sha97,
author={{S}hapiro,{S}tewart}, 
title={Philosophy of Mathematics. Structure and Ontology},
publisher = {Oxford University Press},
address={New York/Oxford} , 
year ={1997},
pages = {x + 279}}

@article{Shoenfieldomega,
author ="Shoenfield,  Joseph R.",
title ="On a restricted $\omega$-rule",
journal= "Bull. Acad. Polon. Sci. Ser. Sci. Math. Astr. Phys.",
year = {1959},
volume ={ 7},
pages ={405-407}}

@article{Speitel2,
    author = "Speitel, Sebastian G.W.",
    title = "Securing Arithmetical Determinacy",
    journal = "Ergo an Open Access Journal of Philosophy, 11:40",
    pages = "1083-1118",
    year = "2024",    
}

@incollection{Tarskilogcon,
author = "Alfred Tarski",
booktitle = "{Logic, Semantics and Metamathematics: Papers from 1936/1983}",
title="On the concept of logical consequence",
publisher = "Hackett",
address ={Indianapolis},
year = "1983",
note = "1st edition 1956",
pages = {409-420}
}

@article{Topey      ,
author = "Brett Topey",
title = "Can we follow the omega rule?",
journal = "The Philosophical Quarterly",
note = {https://doi.org/10.1093/pq/pqaf059},
year = "2025"
}

@book{Warren,
title ={{Shadows of Syntax. Revitalizing Logical and Mathematical Conventionalism}},
author = {Warren, Jared },
publisher = "{Oxford Unversity Press}",
address = {Oxford}, 
year = {2020} }

@book{Whiteheadpr,
 author ={Whitehead, Alfred North}, 
title = "Process and Reality" , 
year ={1929},
publisher ="Harper",
 pages = "xii + 544"}
\bibliographystyle{alpha}

%
%
%
\end{document}